\documentclass[a4paper,12pt]{article}

\usepackage{fancyhdr}
\usepackage{epsfig}
\usepackage{cite}
\usepackage{amsmath}
\usepackage{theorem}
\usepackage{graphicx}
\usepackage{amssymb}
\usepackage{latexsym}
\usepackage{epic}
\usepackage{amsfonts,amssymb,amsopn}
\usepackage[all,2cell,dvips]{xy} \UseAllTwocells \SilentMatrices

\newcommand{\C}{\mathbb{C}}
\newcommand{\Z}{\mathbb{Z}}
\newcommand{\A}{\mathbb{A}}
\newcommand{\PP}{\mathbb{P}}

\newcommand{\Autom}{\mathrm{Aut}}
\newcommand{\Iden}{\mathrm{Id}}

\newcommand{\Ga}{\mathbb{G}_{\mathrm{a}}}
\newcommand{\Gm}{\mathbb{G}_{\mathrm{m}}}
\newcommand{\diag}{\mathrm{diag}}
\newcommand{\qed}{{\ifhmode\unskip\nobreak\fi\quad\ensuremath\square}}

\newcommand{\Csm}{C_{\mathrm{sm}}}
\newcommand{\Gal}{\mathrm{Gal}}
\newcommand{\QK}{{\cal Q}_{K}}
\newcommand{\tQK}{\widetilde{{\cal Q}_K}}

\newcommand{\E}{\mathbf{E}}
\newcommand{\chara}{\mathrm{char}}
\newcommand{\disc}{\mathrm{disc}}

{\theorembodyfont{\slshape}\newtheorem{prop}{{\textbf Proposition}}}
{\theorembodyfont{\slshape}\newtheorem{thm}[prop]{{\textbf Theorem}}}
{\theorembodyfont{\slshape}\newtheorem{cor}[prop]{{\textbf Corollary}}}
{\theorembodyfont{\slshape}\newtheorem{lemma}[prop]{{\textbf Lemma}}}

\title{Quartic equations and 2-division on elliptic curves}

\author{George H.\ Hitching}

\begin{document}

\maketitle

\abstract{Let $K$ be a field of characteristic different from $2$ and $C$ an elliptic curve over $K$ given by a Weierstrass equation. To divide an element of the group $C$ by $2$, one must solve a certain quartic equation. We characterise the quartics arising from this procedure and find how far the quartic determines the curve and the point. We find the quartics coming from 2-division of $2$- and $3$-torsion points, and generalise this correspondence to singular plane cubics. We use these results to study the question of which degree 4 maps of curves can be realised as duplication of a multisection on an elliptic surface.}

\section{Introduction}

Let $K$ be a field of characteristic different from $2$ and $C$ an elliptic curve over $K$. The multiplication-by-$n$ map on $C$ is a morphism of degree $n^2$. If $C$ is given by a Weierstrass equation then this map is given by so-called \textsl{division polynomials} as studied by Cassels \cite{Cas1949}; see also Lang \cite{Lan1978}, Chap.\ 2, and Silverman \cite{Sil1986}, Chap.\ 3.\\
\par
The inverse procedure, \textsl{$n$-division} on $C$, involves solving an equation of degree $n^2$ which is a polynomial in these division polynomials (see Silverman \cite{Sil1986}, exercise 3.7). In the present article, we study the case $n = 2$. We begin with a brief review of the group law on a plane cubic curve in Weierstrass form ($\S$\ref{bkgrnd}), and derive the quartic equation of a 2-division ($\S$\ref{2dqs}). In $\S$\ref{process}, we check which quartics arise from this procedure; precisely, we have:\\
\\
\textbf{Theorem \ref{first}} \textsl{For a quartic $q(x) = x^{4} + d_{3}x^{3} + d_{2}x^{2} + d_{1}x + d_0$, write}
\[ a(q) := 16d_{0} - 4 d_{2}^{2} + d_{2}d_{3}^{2} + 2d_{1}d_{3} \quad \hbox{\textsl{and}} \quad e(q) := 8d_{1} - 4d_{2}d_{3} + d_{3}^{3}. \]
\textsl{Then $q(x)$ arises from the 2-division of a point on some elliptic curve over $K$ if and only if either $e(q)$ is a nonzero square in $K$ or $e(q)$ and $a(q)$ are both zero.}\\
\par
Moreover, given one of these conditions, we show easily how to find the curve(s) and the point(s), which are unique up to a choice of sign if $e(q)$ is nonzero. In particular, in $\S$\ref{special} we determine the quartics arising from 2-divisions of 2- and 3-torsion points, and points on singular plane cubics. In $\S$\ref{ptatinfty}, we ``homogenise'' the picture to include 2-division of the point at infinity. Two points of interest emerge regarding the function $e(q)$: it displays a certain ``discriminant-like'' property, and its vanishing turns out to be closely linked with the geometry of the associated 2-division.
\par
In $\S$\ref{covers} we use these results to study the question of which degree 4 coverings of complex projective smooth curves can be realised as duplication of a multisection on an elliptic surface. In particular (Corollary \ref{Galoiscovers}), we show that any Galois degree 4 cover admits such an interpretation (not necessarily in a unique way) and make some remarks on the existence of 2-torsion sections of the surface.
\par
The last section is independent of the rest: we observe that the duplication map, given with respect to our preferred Weierstrass form, can be naturally expressed with statistical functions.
\par
In the appendix, we prove an easy technical result on Galois degree 4 extensions required in $\S$\ref{covers}.

\paragraph{Acknowledgements:} Part of this work can be found in my M.\ Sc.\ dissertation \cite{Hit2001}. I am very grateful to Miles Reid, who supervised this work at the University of Warwick, for much mathematical and practical assistance and for suggesting the question. I thank Remke Kloosterman for anwering many questions and for advice which shortened some arguments. Thanks also to Klaus Hulek, Kristian Ranestad and Katell Savidan for helpful discussions. I acknowledge gratefully the financial and practical support of the Foreign and Commonwealth Office and British Council, the Deutsche Forschungsgemeinschaft Schwerpunktprogramm ``Globale Methoden in der Komplexen Geometrie'', the University of Warwick and Leibniz Universit\"at Hannover. I also thank the Matematisk Institutt of Universitetet i Oslo for hospitality.

\section{Group theory of plane cubic curves} \label{bkgrnd}

In this section we review some well-known facts about plane cubic curves. Among the many references are Silverman \cite{Sil1986} and Knapp \cite{Kna1992}.\\
\par
Let $K$ be a field and $C \subset \PP^{2}_K$ an irreducible cubic curve with at least one point. The open subset of smooth points of $C$ forms an Abelian group $\Csm$ under an operation which we first describe geometrically. Fix a base point $O \in \Csm$. For $P_1$, $P_{2} \in \Csm$, let $P_{1} * P_2$ be the third point of intersection of the line $L$ joining $P_1$ and $P_2$ with $C$, or the tangent to $C$ at $P_1$ if $P_{1} = P_2$. Denote $OO$ the point $O * O$. Then $P_{1} + P_2$ is defined to be $( P_{1} * P_{2} ) * OO$. It is well known that this defines a group structure on $\Csm$ with respect to which $O$ is the identity. 
\par
Now suppose that $K$ is not of characteristic $2$. Then $C$ can be embedded in $\PP^{2}_K$ as the zero locus of a \textsl{Weierstrass equation}
\[ Y^{2}Z = X^{3} + aX^{2}Z + bXZ^{2} + cZ^3 \]
for some $a$, $b$, $c \in K$. We set $O = (0:1:0)$. This is a flex, so $OO = O$ and $O$ is the single point of intersection of $C$ with the line $\{ Z = 0 \}$. Then the group law is particularly straightforward: $P_{1} + P_2$ is the reflection of $P_{1} * P_2$ in the $X$ axis and the inverse of $(X:Y:Z) \in \Csm$ is $(X:-Y:Z)$.
\par
To begin with, we restrict our attention to the affine subset $C \backslash \{ O \}$ with coordinates $x = X/Z$ and $y = Y/Z$. It will be convenient to write $F(x) = x^3 + ax^{2} + bx + c$; then $C$ is singular if and only if $F$ has a multiple root. A singularity of $C$ is either
\begin{itemize}
\item a node, if $F(x) = (x - s)^{2}(x - t)$ for some $s, t \in \bar{K}$, or
\item a cusp, if $F(x) = (x - s)^3$ for some $s \in \bar{K}$.
\end{itemize}
In either case, the singular point is $(s,0)$. In the nodal case, possibly after taking a quadratic extension of the base field, $\Csm$ is the algebraic group $\Gm$ associated to the multiplicative group $K^*$ and in the cuspidal case it is $\Ga$, that associated to the additive group $K$.

\section{The quartic equation of a 2-division} \label{2dqs}

The \textsl{duplication} map $\Csm \to \Csm$ is given by $P \mapsto P + P = 2_{C}P$. We write this map in coordinates on the complement of $O$ in $\Csm$ (see for example Knapp \cite{Kna1992} or Silverman \cite{Sil1986}. Choose $P = (x_{1}, y_{1})$ in this affine open subset and let $L$ be the tangent to $C$ at $P$. If $y_{1} = 0$ then $F(x_{1}) = y_{1}^{2} = 0$, so $L$ is vertical and $P$ is a point of order 2. Otherwise, $L$ is given by the equation $y = mx + d$, where
\[ m = \frac{\partial F / \partial x}{\partial F / \partial y}(P) = \frac{3x_{1}^{2} + ax_{1} + b}{2y_1} \]
and $d = y_{1} - mx_1$. We find the third point of $L \cap C$. The $x$ coordinates of $L \cap C$ are the solutions of the equation
\[ x^3 + ax^{2} + bx + c - (mx + d)^2 = 0. \]
Since this is monic, the sum of the roots is minus the coefficient of $x^2$, that is, $m^2 - a$. Two of the roots are equal to $x_1$, so
\begin{equation} 2_{C}P = (m^2 - a - 2x_{1},-(mx_{2} + d)). \label{2P} \end{equation}
Since $y_{1}^{2} = x_{1}^3 + ax_{1}^{2} + bx_{1} + c$, we have
\begin{equation} m^{2} = \frac{ ( 3x_{1}^{2} + 2ax_{1} + b )^2}{4( x_{1}^3 + ax_{1}^{2} + bx_{1} + c ) } \label{msq} \end{equation}
Note that if $C$ is a singular irreducible plane cubic then the map $2_C$ is also defined by (\ref{2P}) on the complement of the singular point.\\
\par
We now show how 2-division of a point on the complement of $O$ in $\Csm$ gives rise to a quartic equation. Let $P_{2} = (x_{2}, y_{2})$ be a point of $\Csm \backslash \{ O \}$ and $P_{1} = (x_{1}, y_{1})$ be a point whose double is $P_2$ (of course, we only know that $P_1$ belongs to $C(\bar{K})$). By (\ref{2P}), we have
\[ x_{2} = m^{2} - a - 2x_{1}, \]
so, using (\ref{msq}), we have $F^{\prime}(x_{1})^{2} - 4 F(x_{1})( x_{2} + a + 2x_{1} ) = 0$. Expanding this, we see that $x_1$ satisfies the quartic equation
\begin{equation} x^{4} - 4 x_{2} x^{3} - ( 2b + 4ax_{2} ) x^{2} - ( 8c + 4bx_{2}) x + ( b^{2} - 4ac - 4cx_{2} ) = 0. \label{2dq} \end{equation}
We call this the \textsl{quartic equation of the 2-division} $(C, P_{2})$.\\
\\
\textbf{Remark:} Since the inverse of $(x_{2},y_{2})$ in $\Csm$ is $(x_{2}, -y_{2})$, the quartic equation of $(C, -P_{2})$ is the same as that of $(C, P_{2})$.\\
\par
Next, although the quartic (\ref{2dq}) arose from 2-division of a smooth point of $C$, nothing prevents us from substituting the data of a singular point on a nodal or cuspidal curve into (\ref{2dq}) and obtaining a quartic over $K$. We will examine this and other special cases later.

\section{The 2-division of a generic quartic} \label{process}

Now write $\QK$ for the set of all monic quartics over $K$. This is isomorphic to $\A^{4}_K$ by
\[ x^{4} + d_{3}x^{3} + d_{2}x^{2} + d_{1}x + d_{0} \mapsto ( d_{3}, d_{2}, d_{1}, d_{0} ) . \]
We want to determine which $q \in \QK$ arise from 2-divisions. The coefficients of (\ref{2dq}) are polynomials in $x_2$, $a$, $b$ and $c$. Given $q(x) = x^{4} + d_{3}x^{3} + d_{2}x^{2} + d_{1}x + d_0$, we want to solve the system
\[ d_{3} = - 4 x_{2}, \quad d_{2} = - 2b - 4ax_{2} \]
\[ d_{1} = - 4bx_{2} - 8c, \quad d_{0} = b^{2} - 4ac - 4cx_{2} \]
for $x_2$, $a$, $b$ and $c$ and construct a plane cubic curve $C$ together with point(s) $(x_{2}, \pm y_{2})$ whose 2-division quartic(s) are $q(x)$. Firstly,
\begin{equation} x_{2} = -\frac{d_3}{4}. \label{expx2} \end{equation}
Solving further, we find that
\begin{equation} b = \frac{a d_{3} - d_{2}}{2} \quad \hbox{and} \quad c = \frac{a d_{3}^{2} - (2d_{1} + d_{2} d_{3})}{16} \label{expbc} \end{equation}
and
\begin{equation} \left( 8 d_{1} - 4d_{2}d_{3} + d_{3}^{3} \right) a = 16d_{0} - 4 d_{2}^{2} + d_{2}d_{3}^{2} + 2d_{1}d_{3}. \label{expa} \end{equation}
We define
\[ a(q) := 16d_{0} - 4 d_{2}^{2} + d_{2}d_{3}^{2} + 2d_{1}d_{3} \]
and
\[ e(q) := 8d_{1} - 4d_{2}d_{3} + d_{3}^{3}; \]
the latter quantity has two interesting properties that we will presently see. When either $e(q)$ is not zero or both $a(q)$ and $e(q)$ are zero, we can solve for $a$ (in the latter case there are infinitely many solutions) and write down the curve $C = \{ y^{2} = x^{3} + ax^{2} + bx + c \}$.
\par
In order that $q(x)$ arise as the 2-division quartic of a point on $C_{q}(K)$, there is one further condition: there must be a $K$-rational point of $C$ with $x$ coordinate $x_{2} = -d_{3}/4$. This is equivalent to $F(-d_{3}/4)$ being a square in $K$. This brings us to the first interesting property of $e(q)$: by (\ref{expx2}), (\ref{expbc}) and (\ref{expa}) we find that
\begin{equation} F \left( \frac{-d_3}{4} \right) = \frac{4d_{2}d_{3} - 8d_{1} - d_{3}^{3}}{64} \label{Fx2} \end{equation}
which is none other than $-e(q)/64$. Thus it is necessary that $-e(q)$ be a square in $K$. If $-e(q) = \alpha^2$ for some nonzero $\alpha \in K$ then we get two points $\pm P_{2} = (x_{2}, \pm \alpha/8 )$ whose 2-division quartics are $q(x)$. If $e(q)$ vanishes then the situation is more delicate. We will study this in more detail in the next section.\\
\\
Summing up this section, we have:
\begin{thm} A quartic $q(x) = x^{4} + d_{3}x^{3} + d_{2}x^{2} + d_{1}x + d_0$ arises from the 2-division of a point on some elliptic curve over $K$ if and only if either $e(q)$ is a nonzero square in $K$ or $e(q)$ and $a(q)$ are both zero. \label{first} \end{thm}

\newpage

\noindent \textbf{Remarks}
\begin{enumerate}
\item If $e(q)$ is nonzero then $q$ always represents a 2-division over some quadratic extension of $K$.
\item This analysis could be done somewhat more easily with the assumption that $F(x)$ have no quadratic term (this was the approach taken in \cite{Hit2001}) but this obscures certain points and also requires that $\chara(K)$ be different from $3$. Also, somewhat surprisingly, the calculations with the present approach are only slightly more complicated: there are many agreeable cancellations.
\end{enumerate}

\section{Special 2-division quartics} \label{special}

In this section we study some 2-division quartics with particular properties. We characterise quartics coming from 2-divisions of 2-torsion points, and see how the situation degenerates when the curve may be singular. The classification may be summed up as follows (compare with Silverman \cite{Sil1986}, exercise 3.7):
\begin{thm} Let $q$ be the 2-division quartic of a pair $(C,P_{2})$. Then we have the following classification:
\begin{itemize}
\item $\disc(q) \neq 0$ if and only if $C$ is smooth and $P_2$ is not 2-torsion.
\item $q(x) = (x-s)^{2}(x-t)^2$ for distinct $s,t \in K$ if and only if $P_2$ is a 2-torsion point.
\item $q(x) = (x-s)^{2}(x-t)(x-u)$ for distinct $s, t, u \in \bar(K)$ if and only if $C$ is nodal and $P_2$ is smooth and not 2-torsion.
\item $q(x) = (x-s)^{3}(x-t)$ for distinct $s, t \in K$ if and only if $C$ is cuspidal and $P_2$ is a smooth point.
\item $q(x) = (x-s)^4$ if and only if $P_2$ is a singularity of $C$.
\end{itemize} \label{classification} \end{thm}
\textbf{Proof}\\
We split this into several steps.
\paragraph{Step 1:} Suppose $C$ is smooth. Then we claim that either $q$ has four distinct roots or $q$ has two double roots, and that the latter happens if and only if $P$ is a 2-torsion point.
\par
Suppose two roots of $q$ coincide. Since $C$ is smooth, there must be two distinct halves of $P_2$ with the same $x$ coordinate. These must be of the form $\pm P_1$, so
\begin{equation} P_{2} = 2_{C}P_{1} = 2_{C}(-P_{1}) = - 2_{C}P_{1} = -P_{2}, \label{hottp} \end{equation}
whence $P_2$ is of order 2, and the other two halves of $P_2$ also have the same $x$ coordinate, which is different from that of $P_1$ and $-P_1$. Hence $q$ has two double roots.
\par
Conversely, if $P_2$ is a 2-torsion point then by (\ref{hottp}) its halves occur in pairs which have the same $x$ coordinate. Since $C$ is smooth, there are two such pairs, so $q$ has two double roots.
\paragraph{Step 2:} Suppose $C$ is nodal (resp., cuspidal) with singular point $(s,0)$. Then $(x-s)^2$ (resp., $(x-s)^3$) divides $q$.
\par
If $C$ is nodal then $C$ is given by $y^{2} = F(x) = (x-s)^{2}(x-t)$ where $t \neq s$. We noticed earlier, using (\ref{msq}), that the 2-division quartic of $(C, P_{2})$ is given by
\[ F^{\prime}(x)^{2} - 4 F(x)( x_{2} + a + 2x ) = 0. \]
Since $F^{\prime}(x) = 3(x-s)$(linear factor), clearly $(x-s)^2$ divides $q(x)$. The cuspidal case is similar.
\paragraph{Step 3:} Suppose $C$ is nodal and that $P_2$ is not the singularity. Then the roots of $q(x)$ which do in fact correspond to halves of $P_2$ in $\Csm$ are the roots of
\[ q_{1}(x) := \frac{q(x)}{(x-s)^2}. \]
These are distinct from $s$ since the only point on $C$ with $x$ coordinate $s$ is the node. Moreover, one sees as in Step 1 that they coincide if and only if $P_2$ is the 2-torsion point of $\Csm$.
\paragraph{Step 4:} Suppose $C$ has a cusp at $(s,0)$ and $P_2$ is a smooth point. Step 2 shows that $q$ has a triple root $s$. The single half of $P_2$ in $\Csm$ must have $x$ coordinate different from $s$, so $q(x) = (x-s)^{3}(x-t)$ for distinct $s, t \in K$.
\par
Conversely, if $q(x)$ has a triple root and a simple root then $C$ must be singular by Step 1, and it cannot be nodal by Step 3, so it is cuspidal. 
\paragraph{Step 5:} Lastly, we claim that $q(x) = (x-s)^4$ if and only if $P_2$ is the singular point $(s,0)$.
\par
A tedious but straightforward computation shows that substitution of the data of a singular curve and the singular point $(s,0)$ into (\ref{2dq}) yields the quartic $(x-s)^4$. For the converse, we use a second interesting property of the function $e(q)$, which originally came to light during experiments with Maple. For a quartic $q(x) = (x-s)(x-t)(x-u)(x-v)$, a calculation shows that
\begin{equation} e(q) = -(s+t-u-v)(s-t+u-v)(s-t-u+v). \label{eroots} \end{equation}
From this it follows that $e \left( (x-s)^{4} \right) = 0$. Therefore, for any $a \in K$, the singularity $P_{2} = (s,0)$ of the plane cubic
\[ y^{2} = (x-s)^{2}(x-(a-2s)) = x^{3} + ax^{2} + \cdots \]
has ``2-division quartic'' $(x-s)^4$. We observe that $P_2$ is a cusp if $3s = a$ and a node otherwise.\\
\\
The theorem follows. \qed\\
\par
We now give a result linking the two interesting properties of the function $e(q)$:
\begin{cor} Let $q(x)$ be a 2-division quartic which is not of the form $(x-s)^4$. Then the following are equivalent:
\begin{enumerate}
\renewcommand{\labelenumi}{(\arabic{enumi})}
\item $e(q) = 0$.
\item $q$ is the 2-division quartic of a nontrivial 2-torsion point on a plane cubic which may be singular.
\item $q$ has two double roots.
\end{enumerate}
\label{2tor}
\end{cor}
\textbf{Proof}\\
$(1) \Leftrightarrow (2)$: Let $(C, \pm P_{2})$ be a pair of 2-divisions giving rise to the quartic $q$, where $\pm P_{2} = (x_{2}, \pm y_{2})$. By (\ref{Fx2}), we have
\[ y_{2}^{2} = F(x_{2}) = \frac{-e(q)}{64}, \]
so $e(q) = 0$ if and only if $y_{2} = 0$, equivalently $P_2$ is a nontrivial $2$-torsion point.\\
\\
The equivalence $(2) \Leftrightarrow (3)$ follows from Theorem \ref{classification}. \qed\\
\\
\noindent \textbf{Remarks on the function $e$}
\begin{enumerate}
\item The expression (\ref{eroots}) shows that $e(q(x)) = e(q(x + \alpha))$ for all $\alpha \in K$, so $e$ is invariant under translation of the roots of $q$.
\item We notice from (\ref{eroots}) that $e(q)$ vanishes if and only if the sum of some pair of roots of $q$ is equal to the sum of the other pair; for example, if $q$ has two double roots, or if $q(x)$ is of the form $(x^{2} - s^{2})(x^{2} - t^{2})$. Thus $e$ has a ``discriminant-like'' property.
\item There are certain configurations of points of $\A^{1}_K$ which do not occur as $x$ coordinates of the halves of a point on any elliptic curve in Weierstrass form. For example, the configuration $\pm s, \pm t$ does not occur unless $s = t \neq 0$. Also, if $C$ is nodal and $P_2$ is not of order two, then Corollary \ref{2tor} shows that $e(q)$ is not zero. In particular, we have $2s \neq t + u$ by (\ref{eroots}).\end{enumerate}
\quad \\
We discuss one other special type of 2-division quartic.
\begin{lemma}
Let $q(x) = x^{4} + d_{3}x^{3} + d_{2}x^{2} + d_{1}x + d_0$ be the 2-division quartic of $(C, \pm P_{2})$. Then $P_2$ is a $3$-torsion point if and only if $-d_{3}/4$ is a simple root of $q$.
\end{lemma}
\textbf{Proof}\\
Suppose $P_2$ is a $3$-torsion point of $\Csm$. This means that it is a smooth point of $C$ which is a flex. By definition, $P_{2} * P_{2} = P_2$ and $P_{2} + P_{2} = -P_2$ (whence indeed $3P_{2} = O$). Thus one of the halves of $P_2$ in $\Csm$ is $-P_2$, which has the same $x$ coordinate $x_2$ as $P_2$, so $q(x_{2}) = 0$. But by (\ref{expx2}) we have $x_{2} = -d_{3}/4$. For the rest, we have:
\begin{itemize}
\item $P_2$ is not a $2$-torsion point.
\item $-P_2$ belongs to $\Csm$, so is not a singularity of $C$.
\end{itemize}
Therefore, $-d_{3}/4$ is not a multiple root of $q$ by Theorem \ref{classification}.\\
\par
Conversely, suppose that $-d_{3}/4$ is a simple root of $q$. By Theorem \ref{classification}, it does not correspond to a singularity of $C$. Therefore, since $x_{2} = -d_{3}/4$ by (\ref{expx2}), some half $P_1$ of $P_2$ in $\Csm$ has the same $x$ coordinate as $P_2$. This shows that $P_1 = \pm P_2$, so $2_{C} P_{1} = \pm P_1$. Since $P_1$ belongs to the affine subset $C \backslash \{ O \}$, it is not the identity in $\Csm$, so $2_{C} P_{1} \neq P_1$. Therefore $2_{C} P_{1} = - P_1$, so $P_1$ is a 3-torsion point. Thus $2_{C} P_{1} = P_2$ is also a 3-torsion point. \qed

\section{The point at infinity} \label{ptatinfty}

In this section we extend the correspondence between quartics over $K$ and 2-divisions on plane cubics to include pairs $(C, P_{2})$ where $P_2$ may be the point at infinity $O$. We begin by replacing $\QK$ with the set $\tQK$ of nonzero homogeneous quartics over $K$ with roots in $\PP^1$, up to scalar multiple,
\[ \left\{ \sum_{i=0}^{4} d_{i}X^{i}Z^{4-i} : d_{i} \in K, \hbox{ not all zero } \right\} / K^{*}, \]
a copy of $\PP^{4}_K$ with homogeneous coordinates $(d_{4}:d_{3}:d_{2}:d_{1}:d_{0})$. We consider curves given by equations of the form
\[ Y^{2}Z = X^{3} + aX^{2}Z + bXZ^{2} + cZ^{3} =: F_{h}(X:Z) \]
and write $F^{\prime}_{h}(X:Z) := 3X^{2} + 2aXZ + bZ^{2}$. Now one checks that the double of the point $(X : Y : Z)$ is
\begin{equation} \left( 2 Y_{1} \left( F^{\prime}_{h}(X:Z)^{2} - 4 F_{h}(X:Z)(aZ + 2X) \right) : Y_{2} : 8Y_{1} F_{h}(X_{1}:Z_{1}) Z_{1} \right) \label{dupl} \end{equation}
where $Y_2$ is a function of $(X_{1} : Y_{1} : Z_{1})$. Thus one might like to say that the ``homogenised 2-division quartic'' of the pair $\left( C, ( X_{2} : \pm Y_{2} : Z_{2} ) \right)$ is given by
\[ (X_{2} : Z_{2} ) = \left( F^{\prime}_{h}(X:Z)^{2} - 4 F_{h}(X:Z)(aZ + 2X) : 4 F_{h}(X:Z) Z \right) \]
but this does not make sense for the point $O$, so we need to modify our interpretation.
\par
Let $P_{2} = (x_{2}:y_{2}:1)$ be a point of $C \backslash \{ O \}$ and $q(x)$ its 2-division quartic. We notice that the roots of $q$ correspond to the $x$ coordinates of the halves of $P_2$ under the projection $\pi \colon C \backslash \{ O \} \to \A^{1} \subset \PP^1$ from $O$, given by $(x:y:1) \mapsto (x:1)$.
\par
The map $\pi$ is a priori not defined at $O$, but since $C$ is smooth at $O$ and $\PP^1$ is smooth, it extends to a morphism $C \to \PP^1$, which we also denote $\pi$. Blowing up $\PP^2$ at $O$, one checks easily that $\pi(O) = (1:0)$. Thus it makes sense to define the \textsl{homogenised 2-division quartic of $P_2$} as the quartic with solutions in $\PP^1$ which are the images of the halves of $P_2$ by $\pi$. This quartic is given by
\begin{multline*} ( F^{\prime}_{h}(X:Z)^{2} - 4 F_{h}(X:Z)(aZ + 2X) : 4 F_{h}(X:Z)Z ) = \pi(P_{2}) \\
= \left\{ \begin{array}{l} ( X_{2} : Z_{2} ) \hbox{ if } P_{2} \neq O \\ (1:0) \hbox{ if } P_{2} = O. \end{array} \right. \label{homog2dq1} \end{multline*}
When $Z_{2} \neq 0$ then $Z$ is not zero in any solution of this equation, so we can divide by $Z_{2}Z^4$ and we get the 2-division quartic given earlier (note that $x = X/Z$ and $x_{2}=X_{2}/Z_{2}$). On the other hand, if $P_{2} = O$ then we obtain the quartic
\begin{equation} F_{h}(X:Z)Z = 0 \label{ptatinftyquartic} \end{equation}
whose solutions are $(1:0) \in \PP^1$ and, when $C$ is smooth, the projections of the 2-torsion points of $C$, as they should be. We now say more about the singular case.

\subsection*{Homogenising earlier results}

We would like to extend Theorem \ref{classification} to the completed picture. If $d_4$ is nonzero then we can reduce to the affine situation. Suppose, then, that
\begin{equation} q(X:Z) = d_{3}X^{3}Z  + d_{2}X^{2}Z^{2} + d_{1}XZ^{3} + d_{0}Z^{4} \label{infq} \end{equation}
is a quartic not covered by this patch.
\par
Firstly, the situation described in Theorem \ref{2tor} does not arise here, because (i) the point at infinity is not a nontrivial 2-torsion point, and (ii) if a quartic of the form (\ref{infq}) has two double roots then it is divisible by $Z^2$, so $d_{3} = 0$, but the coefficient of $X^{3}Z$ in the quartic (\ref{ptatinftyquartic}) is nonzero.\\
\\
Suppose, then, that $d_3$ is not zero, so $(1:0)$ is a simple root of $q$.\\
\\
\textbf{Homogenisation of Theorem \ref{classification}} \textsl{A quartic $q(X:Z)$ of the form (\ref{infq}) represents 2-division of the point at infinity on a nodal (resp., cuspidal) plane cubic if and only if $q$ has a double root and a simple root (resp., a triple root) in addition to the root $(1:0)$.}\\
\\
\textbf{Proof}\\
We have $q(X:Z) = Z F_{h}(X:Z)$, so the equation of $C$ can simply be read off. The root $(1:0)$ corresponds to the point at infinity. The other roots correspond to the intersection of $C$ with the $x$ axis; clearly, exactly two (resp., three) of them coincide if and only if the curve is nodal (resp., cuspidal). \qed

\section{Degree 4 coverings of curves} \label{covers}

The preceding sections show that to almost any situation where a quartic equation arises, one can associate a point on an elliptic curve whose geometry reflects some aspects of the situation. An example of such a situation is a degree 4 covering of curves, which is equivalent to a certain quartic field extension. In this section we will use our constructions to study the the following question:\\
\\
Let $\phi \colon X \to Y$ be a degree 4 map of complex projective smooth connected curves. When can $\phi$ be realised as duplication of a multisection of an elliptic surface over $Y$? Precisely, when does there exist an elliptic surface $E \to Y$, a section $\sigma \colon Y \to E$ and a generically injective map $\tau \colon X \to E$ such that the diagram
\begin{equation} \xymatrix{ X \ar[r]^{\tau} \ar[d]_{\phi} & E \ar[d]^{2_{\E}} \\ Y \ar[r]^{\sigma} & E } \label{commdiag} \end{equation}
is commutative?\\
\\
We begin by giving a little background on complex elliptic surfaces. References for this subject include Miranda \cite{Mir1989}, Lecture II and Silverman \cite{Sil1994}, Chapter III.

\subsection{Elliptic surfaces} \label{ellsurfs}
There are a number of different notions of ``elliptic surface over $Y$'' from which we could choose. We will begin with
\begin{enumerate}
\renewcommand{\labelenumi}{(\roman{enumi})}
\item \textsl{smooth minimal elliptic surfaces}: smooth surfaces fibred over $Y$, almost all of whose fibres are elliptic curves and such that no fibre contains curves of self-intersection $-1$, admitting a section not meeting any singular points of singular fibres
\item elliptic curves over the function field $K(Y)$, with a $K(Y)$-rational point
\end{enumerate}
There is a bijection between (isomorphism classes of) the objects (i) and (ii). Later we will consider another model also.

\subsection{Necessary conditions}

Here we will suppose given a commutative diagram of the form (\ref{commdiag}), and deduce some facts about $E$, $\sigma$ and $\phi \colon X \to Y$.

\begin{prop} Suppose we have a commutative diagram of the form (\ref{commdiag}). Then $\sigma$ has no halves in the group of sections of $E$. In particular, it is not the zero section, and so $E$ has a section distinct from the zero section. \label{sigmanotzero} \end{prop}
\textbf{Proof}\\
Suppose we have a section $\sigma^{\prime}$ of $E$ satisfying $2_{E} ( \sigma^{\prime} ) = \sigma$. Then by commutativity of (\ref{commdiag}), the image of $\tau$ would contain the copy $\sigma^{\prime}(Y)$ of $Y$. Since $\tau$ is generically injective, this implies that $\tau(X)$ is not irreducible. But $X$ is smooth and connected, hence irreducible, which is a contradiction.
\par
In particular, since the zero section $\sigma_0$ has a half in the group of sections of $E$, namely itself, $\sigma$ cannot be equal to $\sigma_0$. The proposition follows. \qed\\
\par
Now the fibre of $E$ over the generic point of $Y$ is an elliptic curve $\E$ over the function field $K(Y)$, which has a rational point corresponding to $\sigma_0$. Therefore we can suppose that $\E$ is a plane cubic in $\PP^{2}_{K(Y)}$ described by a Weierstrass equation
\begin{equation} Y^{2}Z = X^{3} + a X^{2}Z + b X Z^{2} + c Z^3 \label{eqnofE} \end{equation}
where $a$, $b$ and $c$ are rational functions on $Y$. By Proposition \ref{sigmanotzero}, the section $\sigma$ defines another $K(Y)$-rational point $(x_{2}, y_{2})$ satisfying the affine equation $y^{2} = x^{3} + a x^{2} + b x + c =: F(x)$.\\
\par
We will work with the model $\bar{E}$ of $\E$ given as the subset
\[ \left\{ \left( (X:Y:Z), P \right) : Y^{2}Z = X^{3} + a(P) X^{2}Z + b(P) X Z^{2} + c(P) Z^{3} \right\} \]
of $\PP^{2}_{\C} \times Y$. This is likely to be singular, but will not be so along any horizontal curves. Therefore, resolving singularities and taking the relatively minimal model yields again the surface $E$. In particular, $\bar{E}$ and $E$ are birational \emph{as elliptic surfaces over $Y$}. This shows that $K(Y)$-rational points of $\E$, sections of $E$ and sections of $\bar{E}$ are all equivalent objects.
\par
Next for any $P \in Y$, we write
\[ F(P)(x) := x^{3} + a(P) x^{2} + b(P) x + c(P) \]
and $F(P)^{\prime}(x)$ for its derivative with respect to $x$.
\par
With respect to the model $\bar{E}$, the map $\tau$ is given by the pair $(\tau_{1}, \phi)$ where $\tau_1$ is a rational function $X \to \PP^2$. This $\tau_1$ can be given generically on the complement of the $Z$ axis by a pair of rational functions $x_{1}, y_{1} \in K(X)$ satisfying 
\[ y_{1}(Q)^{2} = F(\phi(Q))(x_{1}(Q)) \]
where this is defined. Furthermore, by commutativity of (\ref{commdiag}), the function $x_1$ is a root of the quartic equation
\[ \frac{ F(\phi(Q))^{\prime}(x)^2}{4 F(\phi(Q))(x)} - (a + 2x) = x_{2}(\phi(Q)). \]
This is the image of the 2-division quartic of $(\E, \sigma)$ in $\phi^{*}K(Y)[x]$ via the isomorphism $K(Y)[x] \xrightarrow{\sim} \phi^{*}K(Y)[x]$. It has at least one root $x_1$ in the quartic extension $K(X) / \phi^{*}K(Y)$.\\
\\
\textbf{Notation:} For any $s \in \overline{\phi^{*}K(Y)}$, we denote the minimal polynomial of $s$ over $\phi^{*}K(Y)$ by $m_s$.

\begin{lemma} Suppose we have a diagram of the form (\ref{commdiag}). Then one of the following occurs:
\begin{itemize}
\item The section $\sigma$ is not 2-torsion and, via $\phi^{*} \colon K(Y) \hookrightarrow K(X)$, the 2-division quartic of the pair $(\E, \sigma)$ coincides with $m_{x_1}$.
\item The section $\sigma$ is 2-torsion. Then $x_1$ generates a quadratic extension of $\phi^{*}K(Y)$ and $K(X)$ is a quadratic extension of $\phi^{*}K(Y) \left( x_{1} \right)$ generated by $y_1$. In this case the 2-division quartic of $(\E, \sigma)$ coincides with the square of $m_{x_1}$.
\end{itemize}
\label{minpols} \end{lemma}
\textbf{Proof}\\
Suppose $\sigma$ is not a 2-torsion section. We show firstly that $x_1$ generates the extension $K(X) / \phi^{*}K(Y)$. Write $K^{\prime}$ for the subfield of $K(X)$ generated by $x_1$ over $\phi^{*}K(Y)$. We have a sequence of field inclusions
\[ \phi^{*}K(Y) \hookrightarrow K^{\prime} \hookrightarrow K(X), \]
which fix the constant field $\C$. By Silverman \cite{Sil1986}, Theorem II.2.4, there exists a smooth curve $X^{\prime}$, unique up to isomorphism, with function field $K^{\prime}$, and maps $X \to X^{\prime}$ and $X^{\prime} \to Y$ whose composition is equal to $\phi$.
\par
Now we notice that
\[ \frac{- y_{2} - y_1}{x_{2} - x_1} = \frac{F^{\prime}(x_{1})}{2y_1}, \]
since both are the slope of the line in $\PP^{2}_{\overline{K(Y)}}$ passing through $(x_{2}, -y_{2})$ and tangent to $\E$ at $(x_{1}, y_{1})$. Therefore,
\begin{equation} - y_{2} y_{1} - y_{1}^{2} = - y_{2} y_{1} -  F(x_{1}) = \frac{F^{\prime}(x_{1})}{2} \left( x_{2} - x_{1} \right). \label{y1first} \end{equation}
Now we use the fact that $\sigma$ is not 2-torsion. This implies that $y_2$ is not zero in $K(Y)$, so we obtain
\begin{equation} y_{1} = \frac{F^{\prime}(x_{1}) \left( x_{2} - x_{1} \right) + 2 F(x_{1})}{- 2 y_2}. \label{y1} \end{equation}
Now $y_2$ and the coefficients of $F$ all belong to $\phi^{*}K(Y)$, and $x_1$ belongs to $K^{\prime}$. Therefore $y_1$ also belongs to $K^{\prime}$. But this shows that the map $\tau$ factorises via $X^{\prime}$. Since we supposed $\tau$ to be generically injective, we must have $K^{\prime} = K(X)$ and $x_1$ is indeed a primitive element for $K(X)$ over $\phi^{*}K(Y)$.
\par
Next, by Lang \cite{Lan1993}, Proposition V.1.4, the degree of $m_{x_1}$ is equal to $[K(X):\phi^{*}K(Y)] = 4$. Therefore, since the 2-division quartic of $(\E, \sigma)$ is a monic quartic equation over $\phi^{*}K(Y) \cong K(Y)$ which is satisfied by $x_1$, it must be equal to $m_{x_1}$.\\
\par
On the other hand, suppose $\sigma$ is a 2-torsion section of $E \to Y$. In particular, $y_{2} = 0$. By Corollary \ref{2tor}, the 2-division quartic of $(\E, \sigma)$ is of the form $(x-x_{1})^{2}(x-s)^2$ for some $s \in \overline{K(Y)}$. It is easy to check that in fact $(x-x_{1})(x-s)$ has coefficients in $K(Y)$. Now we claim that neither $x_1$ nor $s$ belongs to $K(Y)$. For, by (\ref{y1first}) and since $y_{2} = 0$, we have
\begin{equation} y_{1}^{2} = F(x_{1}) = \frac{F^{\prime}(x_{1})}{2} \left( x_{1} - x_{2} \right). \label{y1again} \end{equation}
If $x_1$ or $s$ were in $K(Y)$ then $y_1$ would belong to a quadratic extension of $K(Y)$, and as above $\tau \colon X \to E$ would factorise via a double cover $X \to X^{\prime}$, contrary to hypothesis. Thus $(x-x_{1})(x-s)$ is irreducible, whence it is clearly the minimal polynomial of $x_1$ over $\phi^{*}K(Y)$.
\par
Next, we notice that if $F(x_{1})$ were a square in $\phi^{*}K(Y) \left( x_{1} \right)$ then, by (\ref{y1again}), both $x_1$ and $y_1$ would belong to a quadratic extension of $\phi^{*}K(Y)$ and again $\tau$ would factorise via a double cover $X \to X^{\prime}$. Hence $F(x_{1})$ is not a square in $\phi^{*}K(Y)$. This shows that $y_1$ is of degree 2 over $\phi^{*}K(Y) \left( x_{1} \right)$. Thus $\phi^{*}K(Y)(x_{1}, y_{1})$ has degree 4 over $\phi^{*}K(Y)$. Since it is contained in $K(X)$, they must be equal. \qed\\
\\
\textbf{Remark:} If $\sigma$ is a 2-torsion section then the image of $\tau$ is $(-1)_E$-invariant. Thus the involution $(-1)_E$ defines an involution of $X$ over $Y$. In fact, on $\tau(X)$, this coincides with translation by the 2-torsion point $\sigma$. In $\S$\ref{Galois}, we will investigate further the relationship between 2-torsion sections of $E$ and automorphisms of the covering $\phi \colon X \to Y$.

\subsection{A sufficient condition}

Let $\phi \colon X \to Y$ be a degree 4 map of complex projective smooth connected curves. Here we show that the vanishing of $e(m_{x_1})$ is a deciding factor in the question of whether a commutative diagram of form (\ref{commdiag}) exists. The issue of whether $e(m_{x_1})$ is a square turns out to be less of an obstacle.

\begin{thm} Let $\phi$ be as above, and suppose there is a primitive element $s$ for the extension $K(X) / \phi^{*}K(Y)$ such that $e(m_{s}) \neq 0$. Then there exist an elliptic surface $\pi \colon E \to Y$ with a section $\sigma \colon Y \to E$, not 2-torsion, and a map $\tau \colon X \to E$ such that the diagram (\ref{commdiag}) commutes. \label{existence} \end{thm}
\textbf{Note:} The surface $E$ will in general be far from unique, even if we require that it be smooth and minimal.\\
\\
\textbf{Proof}\\
Let $s$ be a primitive element for $K(X) / \phi^{*}K(Y)$ such that $e(m_{s})$ is nonzero. If $-e(m_{s}) = \alpha^2$ for some $\alpha \in \phi^{*}K(Y)$ then we set $x_{1} := s$ and $y_{2} := \alpha/8$ (this involves a choice of sign) and, as in $\S$\ref{process}, from the coefficients of $m_s$ we build an irreducible plane cubic curve $\E$ over $K(Y) \cong \phi^{*}K(Y)$ given by the equation $y^{2} = x^{3} + a x^{2} + b x + c =: F(x)$, together with a section $\sigma := (x_{2} , y_{2})$ whose 2-division is represented by $m_s$. Since $\phi$ is separable, the roots of $m_s$ are distinct so in fact $\E$ is smooth by Theorem \ref{classification}. We consider the surface $\bar{E}$ as before. Since by (\ref{Fx2}) we have $F(x_{2}) = -e(m_{s})/64 \neq 0$, the function $y_2$ is nonzero. This fact implies that $\sigma$ is not 2-torsion, and also allows us to construct $y_1$ as in (\ref{y1}). Then the point $(x_{1}, y_{1})$, viewed as an element of
\[ \PP^{2}_{K(X)} \subset \PP^{2}_{\overline{\phi^{*}K(Y)}} \cong \PP^{2}_{\overline{K(Y)}}, \]
satisfies the equation of $\E(\phi^{*}K(Y))$. We define a map $\tau \colon X \to \bar{E}$ by
\[ Q \mapsto \left( (x_{1}(Q): y_{1}(Q) : 1), \phi(Q) \right). \]
Although $\sigma$ and $\tau$ are a priori only rational maps, they extend to morphisms since $X$ and $Y$ are smooth curves over $\C$. For each $Q \in X$ such that $\tau(Q)$ and $\sigma(\phi(Q))$ belong to the complement of the identity in a smooth fibre of $\bar{E}$, we have
\begin{align*} 2_{\bar{E}} ( \tau(Q) ) &= 2_{\bar{E}} \left( ( x_{1}(Q), y_{1}(Q) ), \phi(Q) \right) \\
&= \left( ( \phi^{*}x_{2}(Q) , \phi^{*}y_{2}(Q) ) , \phi(Q) \right) \hbox{ by construction} \\
&= \left( ( x_{2}(\phi(Q)) , y_{2}(\phi(Q)) ) , \phi(Q) \right)\\
&= \sigma(\phi(Q)). \end{align*}
Thus $2_{\bar{E}} \circ \tau = \sigma \circ \phi$ on an open subset of $X$, hence on all of $X$. Resolving singularities of $\bar{E}$ and taking the minimal smooth model $E$, we obtain the diagram (\ref{commdiag}) as required.\\
\par
On the other hand, suppose $-e(m_{s})$ is not a square in $K(Y)$. Write $\varepsilon := e(m_{s})$ for brevity and let $t$, $u$ and $v$ be the other roots of $m_s$. Then one easily checks that the monic quartic
\[ (x - \varepsilon s)(x - \varepsilon t)(x - \varepsilon u)(x - \varepsilon v) \]
is irreducible over $K(Y)$. Hence it is the minimal polynomial of $\varepsilon s$. Moreover, by (\ref{eroots}) we have
\begin{align*} e (m_{\varepsilon s}) &= -( \varepsilon s + \varepsilon t - \varepsilon u - \varepsilon v) ( \varepsilon s - \varepsilon t + \varepsilon u - \varepsilon v) ( \varepsilon s - \varepsilon t - \varepsilon u + \varepsilon v) \\
 &= -\varepsilon^{3} ( s + t - u - v) ( s - t + u - v) ( s - t - u + v) \\
 &= \varepsilon^4 \end{align*}
which is a square in $\phi^{*}K(Y)$. But since $\varepsilon \in \phi^{*}K(Y)$, the product $\varepsilon s$ is still a primitive element for the extension $K(X) / \phi^{*}K(Y)$. We set $x_{1} := \varepsilon s$, and from $m_{\varepsilon s}$ we can construct the elliptic surface over $Y$ and the maps $\sigma$ and $\tau$ as above. \qed

\subsection{Galois degree 4 covers}\label{Galois}

We now say more about the case when $\phi$ is a Galois cover. Write $\Autom(\phi)$ for the group of automorphisms $\alpha$ of $X$ such that $\phi \circ \alpha = \phi$.

\begin{cor} If $\phi$ is a Galois cover then a commutative diagram of the form (\ref{commdiag}) can always be found, such that $\sigma$ is not a 2-torsion section. Such an $E$ always has at least one 2-torsion section, and if $\Autom(\phi) = \Z_4$ then there is exactly one. \label{Galoiscovers} \end{cor}
\textbf{Proof}\\
The first statement will follow from Theorem \ref{existence} and the following technical result:
\begin{prop} Let $K$ be a field of characteristic different from 2, and suppose $L/K$ is a Galois extension of degree 4. If $\Gal(K/L) = \Z_4$, assume that $K$ contains $\sqrt{-1}$. Then there exists at least one primitive element $s$ for the extension whose minimal polynomial $m_s$ over $K$ satisfies $e(m_{s}) \neq 0$. \label{techprimelt} \end{prop}
\textbf{Proof}\\
This proof is straightforward but rather long. We relegate it to the appendix, so as not to interrupt the story.\\
\par
Since a primitive element $s$ exists with the required property, we can construct the surface $E$ and the maps $\sigma$ and $\tau$ by Theorem \ref{existence}.\\
\par
Now for the second statement. We observe that since $\phi$ is a Galois cover, $K(X) / \phi^{*}K(Y)$ is a Galois extension. In particular, $m_{x_1}$ splits, and there exist four distinct maps $\tau \colon X \to E$ making (\ref{commdiag}) commute. We regard these as points on the elliptic curve $\E(K(X)) \supset \E \left( \overline{\phi^{*}K(Y)} \right) \cong \E \left( \overline{K(Y)} \right)$. The difference of any pair of these points is a nontrivial 2-torsion point on $\E$ defined over the field extension $K(X) / \phi^{*}K(Y)$; in fact, all 2-torsion points of $\E$ are defined over this extension. Such a point has coordinates $(t_{1},0)$ for some $t_{1} \in K(X)$. But $F(t_{1}) = 0$, so $m_{t_1}$ divides $F$. Now $m_{t_1}$ has degree dividing 4 since $t_1$ generates a subfield of the degree 4 extension $K(X)/\phi^{*}K(Y)$. If $m_{t_1}$ has degree 1 then $t_1$ belongs to $\phi^{*}K(Y)$ and $(t_{1},0)$ is defined over $\phi^{*}K(Y)$. On the other hand, if $m_{t_1}$ has degree 2 then $F(x) = (x-z)m_{t_1}$ for some $z \in \phi^{*}K(Y)$, so we obtain the 2-torsion point $(z,0)$ defined over $\phi^{*}K(Y)$.
\par
Since the image $\tau(X)$ is invariant under translation by 2-torsion sections in $E$ (corresponding to points of order 2 in $\E(\phi^{*}K(Y)$), any such translation defines an element of $\Autom(\phi)$. This gives a homomorphism $\E(\phi^{*}K(Y))[2] \to \Autom(\phi)$ which is clearly injective. Hence if $\phi$ is a Galois $\Z_4$ cover then the only possibility for $\E(\phi^{*}K(Y))[2]$ is $\Z_2$. \qed\\
\\
\textbf{Question:} Can it happen that $\phi$ is a Galois $\Z_{2} \times \Z_2$ cover but $\E(K(Y))$ has only one 2-torsion point?

\section{Some statistics} \label{stats}

In this short section we notice that in our choice of coordinates, certain statistical quantities naturally occur in connection with the duplication map. Let $C$ be a smooth irreducible plane cubic given by the Weierstrass equation $y^{2} = x^{3} + ax^{2} + bx + c$ on the complement of $O$, and let $P_{2} = (x_{2} , y_{2})$ be a point of $C \backslash \{ O \}$. Write $q(x)$ for its 2-division quartic
\[ x^{4} - 4 x_{2} x^{3} - ( 2b + 4ax_{2} ) x^{2} - ( 4bx_{2} + 8c ) x + ( b^{2} - 4ac - 4cx_{2} ) = 0. \]

Write $\left( x_{1}^{(i)}, y_{1}^{(i)} \right)$ for the halves of $P_2$ in $C(\bar{K})$ (where $i=1, 2, 3, 4$) and denote $m^{(i)}$ the slope of the tangent line to $C$ at $\left( x_{1}^{(i)}, y_{1}^{(i)} \right)$. Firstly, since the coefficient of $x^3$ in $q(x)$ is minus the sum of the roots, $x_2$ is the average of the $x_{1}^{(i)}$.
\par
Next, taking the average of the equations
\[  y_{2} = m^{(i)} \left( x_{2} - x_{1}^{(i)} \right) + y_{1}^{(i)}, \]
we find that $y_{2} = \overline{y_1} + \overline{m} \cdot \overline{x_1} - \overline{ m x_1 }$. But one knows that $\overline{m} \cdot \overline{x_1} - \overline{ m x_1 }$ is the covariance of the variables $x_1$ and $m$, so we have the following pleasant interpretation of the duplication map:
\[ 2_{C}P_{1} = \left( \overline{x_1} , \overline{y_1} + \mathrm{cov}( m , x_{1} ) \right). \]
It would be interesting to investigate whether similar identities hold for the multiplication-by-$n$ maps on $C$ for higher values of $n$.

\appendix

\section{Primitive elements for degree 4 Galois extensions} \label{proofoftechprimelt}

Here we give the proof of Proposition \ref{techprimelt}.\\
\\
\textbf{Proposition \ref{techprimelt}} \textsl{Let $K$ be a field of characteristic different from 2, and $L/K$ a Galois extension of degree 4. If $\Gal(K/L) = \Z_4$, assume that $K$ contains $\sqrt{-1}$. Then there exists at least one primitive element $s$ for the extension whose minimal polynomial $m_s$ over $K$ satisfies $e(m_{s}) \neq 0$.}\\
\\
\textbf{Proof}\\
Let $s$ be any element of $L$. Then the expression
\[ \prod_{g \in \Gal(L/K)}(x - g \cdot s) \]
is a monic quartic over $L$ which is invariant under $\Gal(L/K)$, so in fact belongs to $K[x] \subset L[x]$. It is clearly irreducible over $K$ if $s$ is primitive, so in this case is the minimal polynomial of $s$ over $K$.
\par
Firstly, suppose $\Gal(L/K) = \Z_{2} \times \Z_2$. Then there exist elements $\alpha, \beta \in L$ satisfying $\alpha^{2} \in K$ and $\beta^{2} \in K$ and such that $L = K( \alpha, \beta )$. Then $\{ 1, \alpha, \beta, \alpha\beta \}$ is a basis for $L$ as a vector space over $K$. The Galois group of $L/K$ is given with respect to this basis by the matrices
\[ \left\{ \Iden, \diag(1,-1,1,-1), \diag(1,1,-1,-1), \diag(1,-1,-1,1) \right\} . \]
Denote $s = a + b \alpha + c \beta + d \alpha \beta$ by the 4-tuple $(a,b,c,d)$. Then by (\ref{eroots}), the quantity $e(m_{s})$ is equal to
\begin{multline*} - \left( (a,b,c,d) + (a,-b,c,-d) - (a,b,-c,-d) - (a,-b,-c,d) \right) \\
\cdot \left( (a,b,c,d) - (a,-b,c,-d) + (a,b,-c,-d) - (a,-b,-c,d) \right) \\
\cdot \left( (a,b,c,d) - (a,-b,c,-d) - (a,b,-c,-d) + (a,-b,-c,d) \right), \end{multline*}
which becomes
\[ \left( 0,0,-4c,0 \right) \cdot \left( 0,4b,0,0 \right) \cdot \left( 0,0,0,4d \right) = -64 b c d \alpha^{2} \beta^{2}. \]
Clearly the element $s_{0} = (1,1,1,1) = 1 + \alpha + \beta + \alpha \beta$ is not contained in any degree 2 intermediate field of $L/K$, so this is a primitive element for $L/K$ with $e(m_{s_0}) = -64 \alpha^{2} \beta^{2} \neq 0$.\\
\par
On the other hand, suppose $\Gal(L/K) = \Z_4$. Then by Lang \cite{Lan1993}, Theorem VI.6.2 (i), there is a primitive element $\alpha$ for $L/K$ which satisfies an equation of the form $x^{4} - k = 0$ for some $k \in K$. Then $\{1, \alpha, \alpha^{2} , \alpha^{3} \}$ is a basis for $L$ as a vector space over $K$. Write $i := \sqrt{-1}$. The Galois group is generated by the map defined on basis elements by $\alpha^{j} \mapsto \alpha^{j} i^j$. The group is given with respect to this basis as
\[ \left\{ \Iden, \diag(1, i, -1, -i ) , \diag(1, -1, 1, -1 ) , \diag(1, -i, -1, i ) \right\} . \]
As before, we represent the element $s = a + b \alpha + c \alpha^{2} +d \alpha^3$ by the 4-tuple $(a,b,c,d)$. By (\ref{eroots}) we see that $e(m_{s})$ is equal to
\begin{multline*} - \left( (a,b,c,d) + (a, ib,-c,-id) - (a,-b,c,-d) - (a, -ib,-c,id) \right) \\
\cdot \left( (a,b,c,d) - (a, ib,-c,-id) + (a,-b,c,-d) - (a, -ib,-c,id) \right) \\
\cdot \left( (a,b,c,d) - (a, ib,-c,-id) - (a,-b,c,-d) + (a, -ib,-c,id) \right), \end{multline*}
which becomes
\[ \left( 0,(2 + 2i)b,0,(2 - 2i)d \right) \cdot \left( 0, 0, 4c, 0 \right) \cdot \left( 0, (2 - 2i)b, 0, (2 + 2i)d \right), \]
that is, $32ck(b^{2} + kd^{2})$. Let $l \in K$ be any nonzero element such that $l^{2} \neq -k$. Then the element $s_{0} := 1 + l\alpha + \alpha^{2} + \alpha^3$ does not belong to the one intermediate field of $L/K$, so is a primitive element for $L/K$ with $e(m_{s_0}) = 32k(l^{2} + k) \neq 0$. This completes the proof of the proposition. \qed

\quad\\
\noindent Institut f\"ur Algebraische Geometrie,\\
Leibniz Universit\"at Hannover,\\
Welfengarten 1,\\
30167 Hannover,\\
Germany.\\
\texttt{hitching@math.uni-hannover.de}


\begin{thebibliography}{99}

\bibitem{Cas1949} Cassels, J.\ W.\ S.: \textsl{A note on the division values of $\wp(u)$}. Proc.\ Camb.\ Philos.\ Soc.\ \textbf{45} (1949), pp.\ 167--172.

\bibitem{Cas1991} Cassels, J.\ W.\ S.: \textsl{Lectures on elliptic curves}. London Math.\ Soc.\ Student Texts 24, Cambridge University Press, Cambridge, 1991.

\bibitem {Hit2001} Hitching, G.\ H.: \textsl{Topics in algebraic curves}. M.\ Sc.\ dissertation, University of Warwick, 2001.

\bibitem {HKW1993} Hulek, K.; Kahn, C.; Weintraub, C.: \textsl{Moduli Spaces of Abelian Surfaces: Compactification, Degenerations and Theta Functions}. De Gruyter Expositions in Mathematics 12, Walter de Gruyter, Berlin, 1993.

\bibitem {Kna1992} Knapp, A.\ W.: \textsl{Elliptic curves}. Mathematical Notes 40. Princeton University Press, Princeton, New Jersey, 1992.

\bibitem {Lan1978} Lang, S.: \textsl{Elliptic curves; Diophantine analysis}. Grundlehren der Mathematischen Wissenschaften 231, Springer-Verlag, Berlin--New York, 1978.

\bibitem {Lan1993} Lang, S.: \textsl{Algebra}, third edition. Addison--Wesley, USA, 1993.

\bibitem {Mir1989} Miranda, R.: \textsl{The basic theory of elliptic surfaces}. Dottorato di Ricerca in Matematica, ETS Editrice, Pisa, 1989.

\bibitem {Rei2000} Reid, M.: \textsl{Elliptic Curves}. Lecture notes from course MA426, University of Warwick, 2000.

\bibitem {Sil1986} Silverman, J.\ H.: \textsl{The arithmetic of elliptic curves}. Graduate Texts in Mathematics 106, Springer-Verlag, New York, 1986.

\bibitem {Sil1994} Silverman, J.\ H.: \textsl{Advanced topics in the arithmetic of elliptic curves}. Graduate Texts in Mathematics 151, Springer-Verlag, New York, 1994.

\end{thebibliography}
\end{document}